\documentclass[11pt]{amsart}

\usepackage{amssymb, amsfonts}
\usepackage{graphicx}

\newtheorem{theorem}{Theorem}[section]

\newtheorem{question}{Question}[section]
\theoremstyle{definition}
\newtheorem{definition}{Definition}[section]
\theoremstyle{example}
\newtheorem{example}{Example}[section]
\theoremstyle{remark}

\numberwithin{equation}{section}

\newcommand{\BB}{\mathbb B}
\newcommand{\CC}{\mathbb C}
\newcommand{\RR}{\mathbb R}

\def\RR{{\mathbb R}}

\def\cC{{\mathcal C}}

\def\cE{{\mathcal E}}
\def\cF{{\mathcal F}}
\def\cG{{\mathcal G}}
\def\cH{{\mathcal H}}

\def\cK{{\mathcal K}}

\def\cR{{\mathcal R}}
\def\cS{{\mathcal S}}

\def\Aut{\hbox{\rm Aut}\,}
\def\im{\hbox{\rm Im}\,}
\def\re{\hbox{\rm Re}\,}

\title[Uniform squeezing property and squeezing function]{On the uniform
squeezing property and the squeezing Function}
\author{Kang-Tae Kim and Liyou Zhang}
\address{(Kim) Center for Geometry and its Applications and Department of
Mathematics, POSTECH, Pohang City 790-784 The Republic of Korea}
\email{kimkt@postech.ac.kr}
\address{(Zhang) Department of Mathematics, Capital Normal University, Beijing,
China}
\email{zhangly@mail.cnu.edu.cn}

\begin{document}
\maketitle

\section{Introduction}

In \cite{LSY-1, LSY-2} and \cite{Y}, the concept called
\it holomorphic-homogeneous-regular \rm and equivalently
\it the uniformly-squeezing\rm, respectively, for complex
manifolds has been introduced.  This concept was essential for estimation of
several invariant metrics.  See the above cited papers for details.
\medskip

Let $\Omega$ be a complex manifold of dimension $n$. The {\it squeezing function}
$\sigma_\Omega:\Omega \to \RR$ of $\Omega$ is defined as follows:
for each $p \in \Omega$ let
$$
{\mathcal F}(p,\Omega) := \{f\colon \Omega \to \BB^n, \hbox{ 1-1 holomorphic},
f(p)=0\},
$$
where:
\begin{itemize}
\item  $\BB^n (p;r) = \{ z \in \CC^n \colon \|z-p\|<r \}$, and
\item $\BB^n = \BB^n (0;1)=\BB^n ((0,\ldots,0); 1)$.
\end{itemize}
Then
$$
\sigma_\Omega (p) = \sup  \{r \colon \BB^n (0,r)\subset f(\Omega),
\hbox{ for some }f \in {\mathcal F}(p,\Omega)\}.
$$
Furthermore, the {\it squeezing constant} $\hat \sigma_\Omega$ for $\Omega$ is defined by
$$
\hat\sigma_\Omega := \inf_{p\in\Omega} \sigma_\Omega (p)    .
$$

\begin{definition}[Liu-Sun-Yau \cite{LSY-1, LSY-2}; Yeung \cite{Y}]
\rm A complex manifold $\Omega$ is called {\it holomorphic homogeneous regular}
(HHR), or equivalently {\it uniformly squeezing} (USq), if $\hat\sigma_\Omega>0$.
\end{definition}

Notice that the property HHR (i.e., USq) is preserved by biholomorphisms.  The
squeezing function and squeezing constants are also biholomorphic invariants.

These concepts have been developed in order for the study of completeness and other
geometric properties such as the metric equivalence of the invariant metrics including
Carath\'eodry, Kobayashi-Royden, Teichm\"uller, Bergman, and Kaehler-Einstein
metrics. It is obvious that the examples of HHR/USq manifolds include bounded
homogeneous domains. In case the manifold is biholomorphic to a bounded domain and
the holomorphic automorphism orbits accumulate at every boundary point, such as in the
case of the Bers embedding of the Teichm\"uller space, again USq/HHR property holds. A
bit less obvious example may be the bounded strongly convex domains (as the majority
of them do not possess any holomorphic automorphisms except the identity map),
proved by S.-K. Yeung \cite{Y}.  But there, one of the most standard examples,
such as the bounded convex domains and the bounded strongly pseudoconvex domains
were left untouched.

Indeed the starting point of this article is to show

\begin{theorem} \label{thm-1}
All bounded convex domains in $\CC^n$ ($n\ge 1$) are HHR (i.e., USq).
\end{theorem}

The concept of squeezing function $\sigma_\Omega$ defined above plays an important
role, and moreover it appeals to us that the further investigations on this function should
be worthwhile.  One immediate observation is that if, $\sigma_\Omega (p)=1$ for some
$p \in \Omega$, then $\Omega$ is biholomorphic to the unit open ball (\cite{DGZ}).  In
the light of studies on the asymptotic behavior of several  invariant metrics of the
strongly pseudoconvex domains, perhaps the following question is natural to pose:

\begin{question} \label{q1}
If $\Omega$ is a bounded strongly pseudoconvex domain in $\CC^n$, would
$\displaystyle{\lim_{\Omega\ni q\to p} \sigma_\Omega(q) = 1}$ hold for
every boundary point $p \in \partial\Omega$?
\end{question}

While we do not know the solution at the time of this writing, fortunately, we are able to
present the following result.

\begin{theorem}
\label{thm-2}
If  $\Omega$ is a  bounded domain in $\CC^n$ with a $\cC^2$ strongly convex
boundary, then $\displaystyle{\lim_{\Omega\ni q\to p} \sigma_\Omega (q) = 1}$
for every $p \in \partial\Omega$.
\end{theorem}

The proof-arguments also clarify and simplify some previously-known theorems;
those shall be mentioned in the final section as remarks.
\medskip

\it Acknowlegements. \rm  This research is supported in part by  SRC-GaiA
(Center for Geometry and its Applications), the Grant 2011-0030044 from The Ministry of
Education, and the research of the first named author is also supported in part by National
Research Foundation Grant 2011-0007831, of South Korea.

\section{Bounded convex domains are HHR/USq manifolds}

The aim of this section is to establish Theorem \ref{1-2} stated below.  Not only does this
theorem cover the case left untreated in \cite{Y}, but our method is different.  (See also
\cite{DGZ} on this matter). Our method uses a version of the ``scaling method in several
complex variables'' initiated by S. Pinchuk \cite{Pinchuk}.  In fact, we use the version
presented in \cite{K}, modified for the purpose of studying the asymptotic boundary
behavior of holomorphic invariants.
\medskip

\begin{theorem} \label{1-2}
Every convex Kobayashi hyperbolic domain in ${\mathbb C}^n$ is \break
HHR/USq.
\end{theorem}

Note that all bounded domains are Kobayashi hyperbolic, and every convex Kobayashi
hyperbolic domain is biholomorphic to a bounded domain. But the bounded realization
may not in general be convex.  In that sense this theorem is more general than
Theorem \ref{thm-1}.
\medskip

\noindent\bf Proof. \rm  We proceed in 5 steps.
\medskip

{\bf Step 1. \it Set-up}. Let $\Omega$ be a convex hyperbolic domain in $\mathbb{C}^n$.
Suppose that $\Omega$ is not HHR/USq. Then there exists a sequence $\{q_j\}$ in
$\Omega$ converging to a boundary point, say $q\in \partial\Omega$ such that
$$
\lim_{j \to \infty} S_\Omega (q_j) = 0.
$$
Needless to say, it suffices to show that such a sequence cannot exist.
\medskip

{\bf Step 2. \it The $j$-th orthonormal frame}.
Let $\langle ~, ~ \rangle$ represent the standard Hermitian inner product of $\CC^n$,
and let $\|v\| = \sqrt{\langle v, v \rangle}$. For every $q \in \CC^n$ and a complex linear
subspace $V$ of $\CC^n$, denote by
$$
B^V (q, r) = \{ p \in \CC^n \colon p-q \in V \hbox{ and } \|p-q\|<r\}.
$$
Now let $q \in \Omega$ and define the positive number $\lambda(q, V)$ by
$$
\lambda (q,V) = \max \{r>0 \colon B^V (q, r) \subset \Omega \}.
$$
This number is finite for each $(q,V)$, whenever $\dim V > 0$,
since $\Omega$ is Kobayashi hyperbolic.

Fix the index $j$ momentarily. Then we choose an orthonormal basis for $\CC^n$,
with respect to the standard Hermitian inner product $\langle ~, ~ \rangle$.
First consider
$$
\lambda_j^1 := \lambda (q_j, \CC^n).
$$
Then there exists $q_j^{1*} \in \partial\Omega$ such that $\|q_j^{1*} - q_j\| =
\lambda_j^1$.  Let
$$
e_j^1 = \frac{q_j^{1*} - q_j}{\|q_j^{1*} - q_j\|}.
$$
Then consider the complex span $\hbox{Span}_\CC \{e_j^1\}$, and let $V^1$ be its
orthogonal complement in $\CC^n$. Then take
$$
\lambda_j^2 := \lambda (q_j, V^1)
$$
and $q_j^{2*} \in \partial\Omega$ such that $q_j^{2*}-q_j \in V^1$ and
$\|q_j^{2*} -q_j\|=\lambda_j^2$. Then let
$$
e_j^2 :=  \frac{q_j^{2*} - q_j}{\|q_j^{2*} - q_j\|}.
$$

With $e_j^1, e_j^2, \ldots, e_j^{\ell}$ and $\lambda_j^1, \lambda_j^2, \ldots,
\lambda_j^\ell$ chosen,
the next element $e_j^{\ell+1}$ is selected as follows. Denote by $V^{\ell}$ the
complex orthogonal complement of
$\hbox{Span}_\CC \{e_j^1, e_j^2, \ldots, e_j^\ell \}$.  Then
$$
\lambda_j^{\ell+1} := \lambda (q_j, V^\ell)
$$
and $q_j^{\ell+1*} \in \partial\Omega$ such that $q_j^{\ell+1*}-q_j \in V^\ell$ and
$\|q_j^{\ell+1*} -q_j\|=\lambda_j^{\ell+1}$.  Let
$$
e_j^{\ell+1} :=  \frac{q_j^{\ell+1*} - q_j}{\|q_j^{\ell+1*} - q_j\|}.
$$
By induction, this process yields an orthonormal set $e_j^1, \ldots, e_j^n$ for $\CC^n$
and the positive numbers $\lambda_j^1, \ldots, \lambda_j^n$.
\medskip

{\bf Step 3. \it Stretching complex linear maps}. Let $\hat e^1, \ldots, \hat e^n$ denote
the standard orthonormal basis for $\CC^n$, i.e.,
$$
\hat e^1 = (1,0,\ldots, 0), \hat e^2 = (0,1,0,\ldots, 0), \ldots, \hat e^n = (0, \ldots, 0, 1).
$$
Define the {\it stretching linear map} $L_j: \CC^n \to \CC^n$ by
$$
L_j (z) = \sum_{k=1}^n \frac{\langle z-q_j, e_j^k \rangle}{\lambda_j^k}~ {\hat e}^k
$$
for every $z \in \CC^n$. Note that, for each $j$, $L_j$ maps $\Omega$
biholomorphically onto its image.
\medskip

{\bf Step 4. Supporting hyperplanes.} Notice that
$$
L_j (q_j) = 0 = (0,\ldots,0), L_j (q_j^{1*}) = \hat e^1, \ldots, L_j (q_j^{n*}) = \hat e^n.
$$
We shall consider the supporting hyperplanes, say $\Pi_j^k$ ($k=1,\ldots,n$), of
$L_j(\Omega)$ at points $L_j(q_j^{k*})$, $k=1,\ldots,n$, repectively.
\medskip

{\it Substep 4.1. The supporting hyperplane $\Pi_j^1$}:  Recall that
$L_j (q_j^{1*}) = \hat e^1 =(1,0,\ldots, 0)$. Due
to the choice of $q_j^{1*}$ the supporting hyperplane of $\Omega$ at $q_j^{1*}$
must also support the sphere tangent to the boundary $\partial\Omega$.
Consequently the supporting hyperplane $\Pi_j^1$ of
$L_j(\Omega)$ must support a smooth surface (an ellipsoid) tangent to
$L_j(\partial\Omega)$ at $\hat e^1$.  Thus the equation for this hyperplane $\Pi_j^1$ is
$$
\re (z_1-1) = 0
$$
(independently of $j$, being perpendicular to $\hat e^1$ consequently).
We also note that
$$
L_j (\Omega) \subset \{ (z_1, \ldots, z_n) \in \CC^n \colon \re z_1 < 1\}.
$$
\medskip

{\it Substep 4.2. The rest of supporting hyperplanes $\Pi_j^k$, for $k\ge 2$}: First
consider the case $k=2$. Then the supporting hyperplane $\Pi_j^2$
passes through $L_j (q_j^{2*}) = \hat e^2 =(0,1,\ldots, 0)$. Since the restriction of
$\Omega$ to $V^1$ contains the sphere in $V^1$ tangent to the restriction of
$\partial\Omega$ at the point $\hat e^2$, the supporting hyperplane $\Pi_j^2$ restricted
to $L_j(V^1)$ takes the equation $\{(z_2,\ldots,z_n) \in \CC^{n-1}\colon
\re (z_2 -1) = 0\}$.
Hence
$$
\Pi_j^2 = \{(z_1, \ldots, z_n)\in\CC^n \colon \re (a_j^{2,1} z_1 + a_j^{2,2} (z_2 -1))
= 0\}
$$
for some $(a_j^{2,1}, a_j^{2,1}) \in \CC^2$ with $\Big|a_j^{2,1}\Big|^2 +
\Big|a_j^{2,2}\Big|^2=1$ and $a_j^{2,2} >0 $.
We also have that
$$
L_j(\Omega) \subset \{ (z_1,  \ldots, z_n) \in \CC^n \colon \re (a_j^1 z_1 +
a_j^2 (z_2 -1))<0\}.
$$
For $k  \in \{3, \ldots, n\}$, one deduces inductively that the supporting hyperplane
$\Pi_j^k$
passes through the point $\hat e^k$, and that
\begin{multline*}
\Pi_j^k = \{(z_1, \ldots, z_n)\in\CC^n \colon \\ \re (a_j^{k,1} z_1 + \cdots +
a_j^{k,k-1}
z_{k-1} + a_j^{k,k} (z_k  -1))=0,
\end{multline*}
with $a_j^{k,k} >0$ and $\sum_{\ell=1}^k \Big|a_j^{k,\ell}\Big|^2=1$.
Also,
\begin{multline*}
L_j (\Omega) \subset   \{(z_1, \ldots, z_n)\in\CC^n \colon \\
\re (a_j^{k,1} z_1 + \cdots + a_j^{k,k-1} z_{k-1} + a_j^{k,k} (z_k  -1))<0\}.
 \end{multline*}
\medskip

{\it Substep 4.3. Polygonal envelopes}:  We add this small substep for convenience.
From the discussion by far in this Step, we have the $j$-th polygonal envelope (of
$L_j(\Omega)$)
\begin{eqnarray*}
\Sigma_j & := &  \{(z_1, \ldots, z_n) \in \CC^n :  \\
& &  \qquad  \qquad \quad  \re z_1 < 1 \\
& &  \qquad \qquad \re (a_j^{2,1} z_1 + a_j^{2,2} (z_2 -1))<0 \\
& &  \qquad \qquad  \qquad  \qquad \vdots \\
& &  \qquad \quad \re (a_j^{n,1} z_1 + \cdots + a_j^{n,n-1} z_{n-1} + a_j^{n,n}
(z_n  -1))<0\}
\end{eqnarray*}
\medskip

{\bf Step 5. Bounded realization}.
Notice that, for every $k \in \{1,\ldots,n\}$,  the disc
$$
D_j^k := \{z=(z_1, \ldots, z_n) \in \CC^n \colon \langle z-q_j, e_j^\ell \rangle = 0,
\forall \ell \neq k; \|z-q_j\|<\lambda_j^k \}
$$
is contained in $\Omega$.  Hence, every $L_j (\Omega)$ contains the discs $D^k := \{
\zeta \hat e^k \colon \zeta \in \CC, |\zeta|<1\}$ for every $k = 1,\ldots, n$.  Since
$\Omega$ is convex and since $L_j$ is linear, $L_j (\Omega)$ is also convex.
Therefore, the ``unit acorn''
$$
A := \{ (z_1, \ldots, z_n) \in \CC^n \colon |z_1|+ \cdots + |z_n| < 1 \}
$$
is contained in $L_j (\Omega)$.  This restricts the unit normal vectors \break $n_j^k :=
(a_j^{k,1}, \ldots, a_j^{k,k}, 0, \ldots, 0) \in \CC^n$ for every $k=2,\ldots, n$.  Namely,
there is a positive constant $\delta>0$ independent of $j$ and $k$ such that $a_j^{k,k}
\ge \delta$ for every $j, k$.

Now taking a subsequence (of $q_j$), we may assume that the sequence of unit vectors
$\{n_j^k\}_{j=1}^\infty$ converges for every $k\in \{2,\ldots,n\}$.  Let us write
$$
\lim_{j\to\infty} n_j^k = n^k = (a^{k,1}, \ldots, a^{k,k}, 0,\ldots, 0)
$$
for each $k = 1,2,\ldots, n$.

Consider the maps
$$
B_j (z_1, \ldots, z_n) = (\zeta_1, \ldots, \zeta_n)
$$
defined by
\begin{eqnarray*}
\zeta_1 & = & z_1, \\
\zeta_2 & = & a_j^{2,1} z_1 + a_j^{2,2} z_2, \\
& \vdots & \\
\zeta_n & = & a_j^{n,1} z_1 + \ldots + a_j^{n,n} z_n. \\
\end{eqnarray*}
Then it follows that
\begin{eqnarray*}
B_j \circ L_j (\Omega) & \subset & B_j (\Sigma_j) \\
& = &
\{(\zeta_1, \ldots, \zeta_n) \in \CC^n \colon \re \zeta_1<1, \re\zeta_2 < a_j^{2,2}, \ldots,
\re\zeta_n < a_j^{n,n}\}
\end{eqnarray*}

Now we consider the Cayley transformation, for each $j$,
$$
\Phi_j (z_1, \ldots, z_n) = \Big( \frac{z_1}{2-z_1}, \frac{z_2}{2a_j^{2,2} - z_2},
\ldots, \frac{z_n}{2a_j^{n,n} - z_n} \Big).
$$
Then $\Phi_j \circ B_j (\Sigma_j) \subset D^n$,
where $D^n$ denote the unit polydisc in $\CC^n$ centered at the origin.   Also, there
exists a positive constant $\delta' \in (0,\delta)$ such that $\Phi_j \circ B_j (\Sigma_j)
\subset D^n$ contains the ball of radius $\delta'$ centered at the origin $0$.

Since $\Phi_j \circ B_j \circ L_j (q_j) = (0,\ldots,0)$ for every $j$,  we now conclude
that the squeezing function satisfies
$$
\sigma_\Omega (q_j) \ge \frac{\delta'}{\sqrt{n}}.
$$
This estimate, which holds for every sequence $q_j$ approaching the boundary, yields
the desired contradiction at last.  Thus the proof is complete. \hfill $\Box$

\section{Boundary behavior of squeezing function on strongly convex domains}

Consider first the following

\begin{definition}
Let $\Omega$ be a domain in $\CC^n$.  A boundary point $p \in \partial\Omega$ is said
to be {\it spherically-extreme} if
\begin{itemize}
\item[(1)] the boundary $\partial\Omega$ is $\cC^2$ smooth in an open neighborhood
of $p$, and
\item[(2)] there exists a ball $\BB^n (c(p);R)$ in $\CC^n$ of some radius $R$, say, centered
at some point $c(p)$ such that $\Omega \subset \BB^n(c(p);R)$ and $p \in \partial\Omega \cap
\partial \BB^n(c(p);R)$.
\end{itemize}
\end{definition}

The main goal of this section is to establish

\begin{theorem}
\label{thm-2g}
If a domain $\Omega$ in $\CC^n$ admits a
spherically-extreme boundary point $p$, say, in a neighborhood of which the
boundary $\partial\Omega$ is $\cC^2$ smooth, then
$$
\lim_{\Omega\ni q \to p} \sigma_\Omega (q)=1.
$$
\end{theorem}

\bf Proof. \rm
Since every boundary point of a $\cC^2$ strongly convex bounded domain is
spherically-extreme, this theorem implies Theorem \ref{thm-2}. The rest of this section is devoted
to the proof of Theorem \ref{thm-2g}, which we shall proceed in seven steps.
\medskip

{\bf Step 1: Sphere Envelopes.}
Let $\Omega$ be a bounded domain in $\CC^n$ with a boundary point
$p \in \partial\Omega$ such that
\begin{itemize}
\item[(\romannumeral 1)]$\partial\Omega \cap B^n (p; r_0)$ is $\cC^2$-smooth
for some $r_0>0$, and
\item[(\romannumeral 2)] $p$ is a spherically-extreme boundary point of $\Omega$.
\end{itemize}

Then there exist positive constants $r_1, r_2$ and $R$ with $r_0>r_1>r_2$ such that
every $q \in \Omega \cap \BB^n (p; r_2)$ admits  points
$b(q) \in \partial\Omega \cap \BB^n (p;r_1)$ and   $c(q) \in \CC^n$
satisfying the conditions
\begin{itemize}
\item[(\romannumeral 3)] $\|q-b(q)\|<\|q-z\|$ for any
$z \in \partial\Omega - \{b(q)\}$, and
\item[(\romannumeral 4)] $\|c(q)-b(q)\|=R$ and $\Omega \subset \BB^n (c(q);R)$.
\end{itemize}

\begin{figure}[h]
\centering
\includegraphics[height=1.9in, width=2.00in]{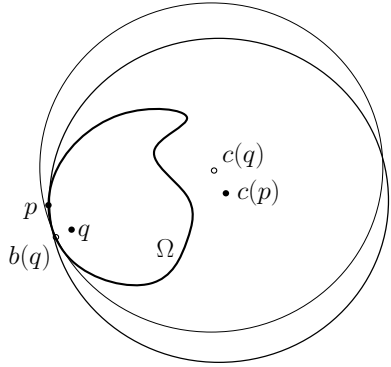}
\caption{\sf Sphere envelopes}
\end{figure}

Notice that (\romannumeral 3) says that $b(q)$ is the unique boundary point that is the closest to $q$, and that the constant $R$ in (\romannumeral 4) is independent of the choice of $q \in \BB^n(p;r_2)$.
\bigskip

{\bf Step 2:  Centering.}
From this stage we shall exploit the familiar notation
\begin{eqnarray} \label{T1}
z  &=& (z_1,\ldots, z_n), \nonumber \\
z' &=& (z_2, \ldots, z_n),  \\
u &=& \re z_1, \nonumber\\
v &=& \im z_1.\nonumber
\end{eqnarray}
\smallskip
For each $q \in \Omega \cap B^n (p, r_2)$, choose a unitary transform $U_q$
of $\CC^n$ such that the map $A_q (z) := U_q (z-b(q))$
satisfies the following conditions:
\begin{equation} \label{T2}
A_q (q) = (\lambda_q, 0, \ldots, 0)
\end{equation}
for some $\lambda_q > 0$, and
\begin{equation} \label{T3}
A_q (\Omega) \subset \BB^n ((R,0,\ldots,0); R)
= \{ z \in \CC^n \colon |z_1-R|^2 + \|z'\|^2 < R^2 \}.
\end{equation}

\begin{figure}[h]
\centering
\includegraphics[height=1.3in, width=3.70in]{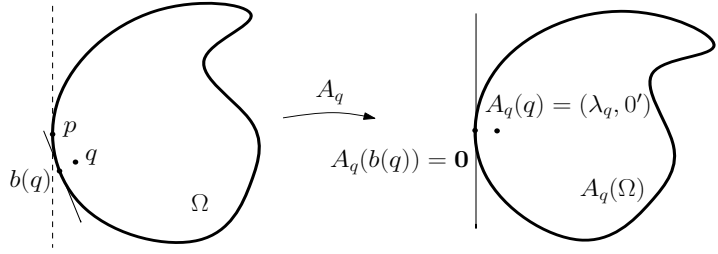}
\caption{\sf The Centering Process}
\end{figure}

Then there exists a positive constant $r_3<r_2$ such that
\begin{multline} \label{T4}
z \in A_q (\Omega)  \cap B^n (0, r_3)
\\ \Leftrightarrow  \|z\| < r_3 ~\hbox{ and }~
2 u > H_{b(q)} (z') + \cK_{b(q)} (v, z') + \cR_{b(q)} (v, z')
\end{multline}
where:
\begin{itemize}
\item $H_{b(q)}$ is a quadratic positive-definite Hermitian form such that there
exists a constant $c_0 > 0$, independent of $q$, satisfying
\begin{equation} \label{T5}
H_{b(q)} (z') \ge c_0 \|z'\|^2
\end{equation}
and
\item there exists a constant $C > 0$, independent of $q \in \BB^n (p; r_3) \cap \Omega$, such that
\begin{equation} \label{T6}
|\cK_{b(q)} (v, z')| \le  C(|v|^2 +|v|\|z'\|),
\end{equation}
whenever $z \in \BB^n (0, r_3)$.  Furthermore, we have
$$
|\cR_{b(q)} (v, z')| = o(|v|^2 + \|z'\|^2).
$$
In particular, the choice of $r_3$ can allow us the estimate
$$
|\cR_{b(q)} (v, z')| \le \frac{c_0}{2}(|v|^2 + \|z'\|^2).
$$

\end{itemize}
Notice that
$$
\lim_{\Omega \ni q \to p} b(q) = p,
\qquad
\lim_{\Omega \ni q \to p} H_{b(q)}(z') = H_p(z'),
$$
and
$$
\lim_{\Omega \ni q \to p} A_q = I \hbox{ (the identity map)}.
$$
This last and an inductive construction yield that for each integer $m>2$ there exists a strictly-increasing integer-valued function $k(m)$  such that
\begin{equation} \label{T7}
\BB^n (0; r_3/(2k(m))) \subset  A_q\big(\BB^n (p; r_3/k(m))\big) \subset \BB^n (0; r_3/m),
\end{equation}
whenever $q \in \BB^n (p, \frac{r_3}{2k(m)})$.
\bigskip

\bf Step 3: The Cayley transform. \rm
The {\it Cayley transform} considered here is the map
\begin{equation} \label{T8}
\kappa (z) := \Big( \frac{1-z_1}{1+z_1}, \frac{\sqrt{2} z_2}{1+z_1}, \ldots, \frac{\sqrt{2} z_n}{1+z_1}\Big),
\end{equation}
well-defined except at points of $Z= \{z \in \CC^n \colon z_1=-1\}$.
Notice that this transform maps the open unit ball $\BB^n (0;1)$ biholomorphically onto the
Siegel half space
\begin{equation} \label{T9}
\cS_0 := \{z \in \CC^n \colon 2 \re z_1 > \|z'\|^2 \}.
\end{equation}

Moreover, $\kappa \circ \kappa = 1$ and consequently, $\kappa (\cS_0) = B^n (0,1)$.  Notice also that, if we denote by ${\bf 1}=(1,0,\ldots)$ and $-{\bf 1}=(-1,0,\ldots)$, then we have
$\kappa({\bf 1}) =(0,\ldots,0)$, $\kappa((0,\ldots,0))={\bf 1}$,
$\kappa(-{\bf 1}) = \infty$  and  $\kappa(\infty) = -{\bf 1}$.
\bigskip

\bf Step 4: Stretching. \rm
Let $q \in \Omega \cap \BB^n (p; \frac{r_3}{2k(m)})$.  If we let $m$ tend to infinity.  Then of course $A_q (q) = (\lambda_q, 0, \ldots, 0)$ approaches $A_q(b(q))=(0,\ldots,0)$ and so $\lambda_q$ approaches zero.  For simplicity, denote by $\lambda = \lambda_q$, suppressing the notation $q$.  But $\lambda$ is still dependent upon $q$.  Note that
\begin{equation} \label{T10}
A_q (\BB^n(c(q); R)) = \{z \in \CC^n \colon 2 R \ \re z_1 > \|z\|^2 \}.
\end{equation}
Define the map $\Lambda_\lambda\colon \CC^n \to \CC^n$ by
\begin{equation} \label{T11}
\Lambda_\lambda (z) := \Big( \frac{z_1}{\lambda}, \frac{z_2}{\sqrt{\lambda}}, \cdots,
\frac{z_n}{\sqrt{\lambda}} \Big),
\end{equation}
the {\it stretching map}, introduced originally by Pinchuk (cf.\ \cite{Pinchuk}).

Recall (\ref{T6}). This stretching map transforms $A_q(\Omega)\cap \BB^n (0; \frac{r_3}{3})$ to the domain
$
\Lambda_\lambda \big(A_q(\Omega) \cap \BB^n (0;\frac{r_3}{3})\big)
$
so that
\begin{eqnarray}
& z &   \in  \Lambda_\lambda \circ A_q(\Omega)
\cap \BB^n \Big(0;\frac{r_3}{\sqrt{\lambda}k(3)}\Big) \label{T12} \\
&  & \Leftrightarrow  \|z\|<\frac{r_3}{\sqrt{\lambda}k(3)} \text{ and } \nonumber\\
&  &\qquad  2u > H_{b(q)}(z') + \frac1\lambda K_{b(q)} (\lambda v, \sqrt{\lambda}z') +
\frac1\lambda \cR_{b(q)} (\lambda v, \sqrt{\lambda}z'). \nonumber
\end{eqnarray}
On the other hand, notice that
$$
\Big\|\frac1\lambda K_{b(q)} (\lambda v, \sqrt{\lambda}z')\Big\|
\le
C\sqrt{\lambda}(\sqrt{\lambda} |v|^2 + |v|\|z'\|)
$$
and that
$$
\Big\|\frac1\lambda \cR_{b(q)} (\lambda v, \sqrt{\lambda}z')\Big\|
\le \frac1\lambda o ((|\lambda v|^2 + \|\sqrt{\lambda}z'\|^2)) = \frac1\lambda o(\lambda)
$$
on $\BB^n(0;\rho)$ for any fixed constant $\rho>0$.  Notice that both terms
approach zero as $\lambda$ tends to zero.  Thus, these terms can become sufficiently
small if we limit $q$ to be contained in $\BB^n (p; \frac{r_3}{2k(m)})$ for some
sufficiently large $m$.
\medskip

\bf Step 5: Set-convergence. \rm
This step is in part heuristic; and the heuristics appearing, especially which concern
set-convergences, in this step are not used in the proof, strictly speaking.
We include this step because they seem to help us to grasp the logical structure of the
proof.  On the other hand, the constructions in (\ref{T13})--(\ref{T15}) shall be used in
the proof-arguments, especially in Step 7.

The main role of the stretching map $\Lambda_\lambda$, as $\lambda \searrow 0$ is to
rescale the domains successively, letting them to converge to the set-limits.

For instance if one considers
$$
\Lambda_\lambda ( A_q (\Omega) \cap B^n (0,r_3))
$$
then, one can see that $\Lambda_\lambda (B^n (0, r_3))$ contains $B^n (0, r_2/\sqrt{\lambda})$,
a very large ball, which exhausts $\CC^n$ successively as $\lambda$ approaches zero.
In the mean time within that large ball, $\Lambda_\lambda (A_q (\Omega))$ is restricted only by
the inequality
$$
2 u > H_{b(q)} (z') + \tilde K_\lambda (v,z')
$$
where $\tilde K_\lambda = o(\lambda)$ is small enough to be negligeable.  One can
imagine that indeed the ``limit domain'' of this procedure should be
\begin{equation} \label{T13}
\widehat\Omega := \{z \in \CC^n \colon 2u > H_p (z')\}.
\end{equation}
Here, of course, $H_p(z')$ is the quadratic positive-definite Hermitian form which appears in the defining inequality of $\Omega$ about the boundary point $p$ (understood as the origin):
$$
2 \re z_1 > H_p (z') + o(|\im z_1|+ \|z'\|^2).
$$
Notice that
$$
\kappa(\widehat\Omega) = \{z \in \CC^n \colon |z_1|^2 + H_p (z') < 1\},
$$
and hence there is a $\CC$-linear isomorphism
\begin{equation}\label{T14}
L\colon \CC^n\to \CC^n
\end{equation}
that maps $\kappa(\widehat\Omega)$ biholomorphically onto the unit ball $\BB^n(0;1)$ with $L({\bf 1})= {\bf 1}$.

Before leaving this step we remark that, since $\Omega \subset \BB^n(c(q);R)$ whenever
$q \in \BB^n (p;r_2)$, $A_q (\Omega) \subset A_q (\BB^n(c(q);R))
= \BB^n ((R,0,\ldots,0);R)$.  This in turn implies that
\begin{eqnarray}
 \Lambda_\lambda \circ A_q(\Omega)
& \subset   &  \Lambda_\lambda\big(\BB^n ((R,0,\ldots,0);R)\big)
\label{T15} \\
& \subset & \cE := \{z\in\CC^n \colon 2R~ \re z_1 > \|z'\|^2\}. \nonumber
\end{eqnarray}
The last inclusion follows by (\ref{T10}).
\bigskip

\bf Step 6: Auxiliary domains. \rm
Let $\delta > 0$ be given.  Consider the domains
\begin{equation}\label{T16}
\cG_\delta := \{z \in \CC^n \colon 2 u > -\delta |v| + (1-\delta) H_{b(q)} (z') \},
\end{equation}
\begin{equation}\label{T17}
\cF_\delta  :=  \{z \in \CC^n \colon 2 u > \delta |v| + (1+\delta) H_{b(q)} (z') \}
\end{equation}
and
\begin{equation}\label{T18}
 \cH_q := \{z\in\CC^n \colon 2 u > H_{b(q)}(z') \},
\end{equation}
in addition to $\widehat\Omega$ and $\cE$ introduced in (\ref{T13}) and (\ref{T15}).
\begin{figure}[h]
\centering
\includegraphics[height=1.8in, width=2in]{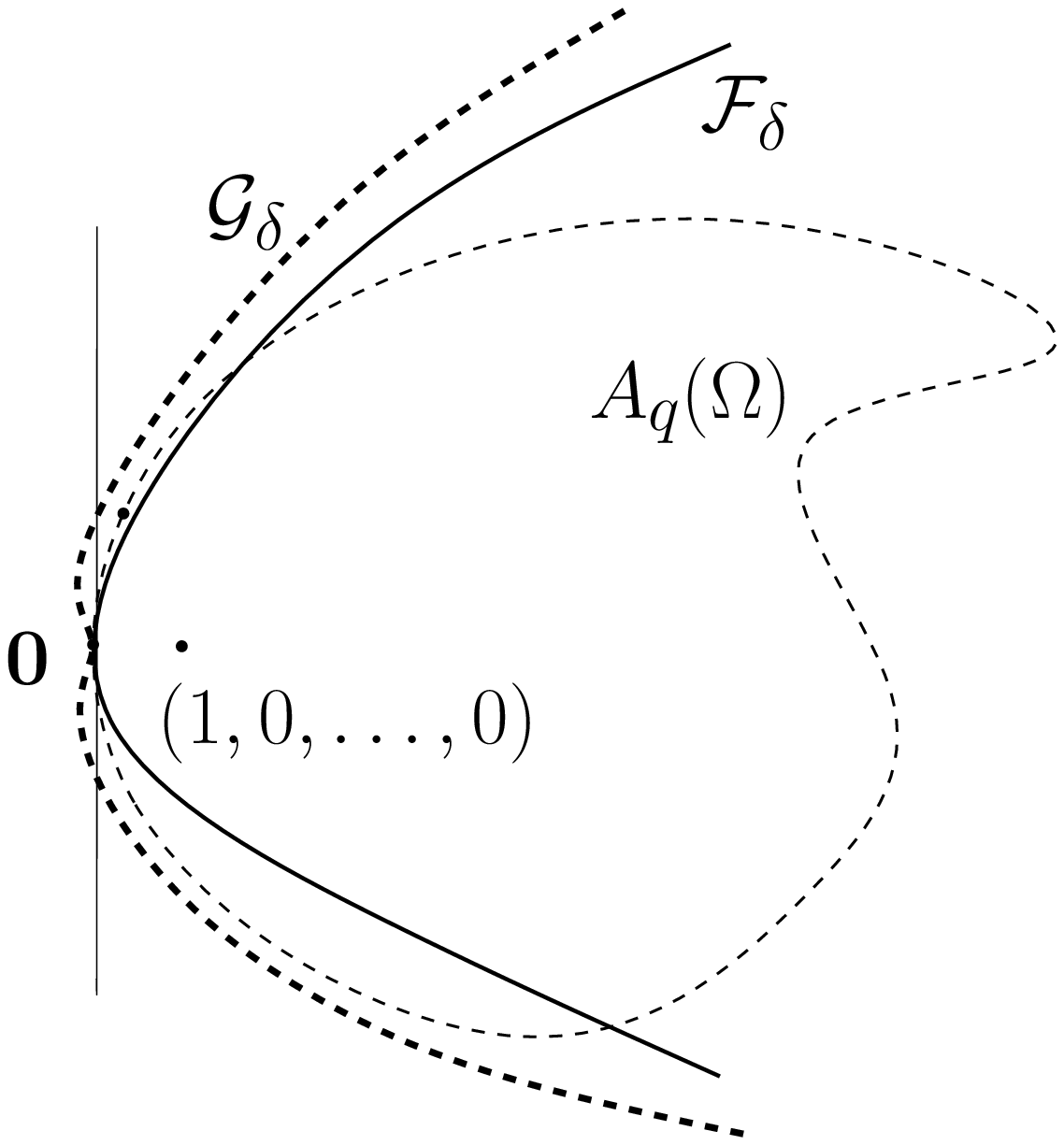}
\caption{\sf Auxiliary domains $\cG_\delta$ and $\cF_\delta$}
\end{figure}
\smallskip

\noindent
A straightforward computation checks that the image $\kappa (\cG_\delta)$ of $\cG_\delta$ via the Cayley transform $\kappa$ introduced earlier is
\begin{equation} \label{T19}
\kappa(\cG_\delta) = \{z \in \CC^n \colon |z_1|^2 - \frac{\delta}2 |z_1-\bar z_1| + (1-\delta)H_{b(q)} (z') < 1\}.
\end{equation}
Hence, there exists $\delta_0 > 0$ that, for every $\delta$ with $0<\delta<\delta_0$,
$\kappa(\cG_\delta)$ is a bounded domain.  Notice also that this domain is arbitrarily
close to the domain $\kappa (\cH_{b(q)})$ as $\delta_0$ becomes arbitrarily small.
It follows therefore that, for every $\epsilon > 0$, there exists $\delta_0>0$ such that
\begin{equation} \label{T20}
L\circ \kappa(\cG_\delta)  \subset \BB^n (0;1+\epsilon)
\end{equation}
\begin{figure}[h]
\centering
\includegraphics[height=1.40in, width=2.90in]{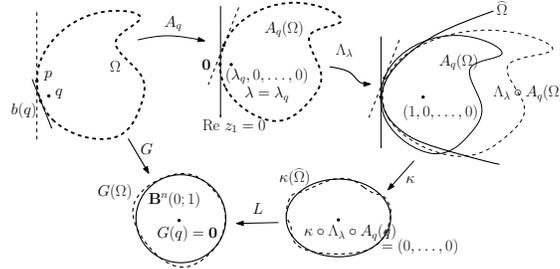}
\caption{$G(\Omega)=L\circ \kappa \circ \Lambda_\lambda \circ A_q(\Omega)$ for $q \sim p$}
\end{figure}
whenever $0<\delta<\delta_0$.  Moreover, observe that the stretching map
$\Lambda_\lambda$ preserves all such domains as
$$
\cF_\delta, \cG_\delta, \widehat\Omega, \cE \text{ and } \cH_q.
$$
Let us now define the expression
\begin{equation} \label{T21}
G(z) := L \circ \kappa \circ \Lambda_\lambda \circ A_q (z)
\end{equation}
for $z \in \CC^n - (\Lambda_\lambda \circ A_q)^{-1} (Z)$. [ The set $Z$ has
been defined in (\ref{T8}). Notice that this expression $G$ depends upon
$q \in \BB^n (0; r_2)$, for instance; see Figure 3 in Step 4 for an illustration.]
In particular, this $G$ maps $\Omega$ onto its image $G(\Omega)$ biholomorphically.

\bigskip

\bf Step 7: Proof of Theorem \ref{thm-2g}. \rm
Our present goal is to show the following
\medskip

\noindent
\bf Claim. \it For any $\epsilon$ with $0<\epsilon<1/2$, there exists an integer $m >0$ such that
\begin{equation} \label{T22}
\BB^n (0; 1-\epsilon) \subset G(\Omega) \subset \BB^n (0; 1+\epsilon)
\end{equation}
whenever $q \in \Omega \cap B^n (p, \frac{r_3}{2k(m)})$. \rm
\medskip

Since $G(q)=0$, this implies that the squeezing function $\sigma_\Omega$ satisfies
$$
\sigma_\Omega (q) \ge \frac{1-\epsilon}{1+\epsilon}.
$$
Notice that this completes the proof of Theorem \ref{thm-2g}.
\medskip

Therefore we are only to establish this claim.
\medskip
Start with $\BB^n(0; 1-\epsilon)$.  Notice first, by the definition of $\cF_\delta$, that
for every $\delta>0$ there exists $m_1>0$ such that
$$
\cF_\delta \cap \BB^n (0;r_2/m) \subset A_q(\Omega) \cap \BB^n (0;r_2/m),
$$
for any $m>m_1$.

Also,
$$
\kappa^{-1}\circ L^{-1}( \BB^n(0; 1-\epsilon)) \subset\subset \kappa^{-1}\circ L^{-1}( \BB^n(0; 1)) = \widehat\Omega.
$$
As discussed in ({T4})--(\ref{T7}), $L \circ \kappa (\cH_q)$ is sufficiently close to
$L\circ \kappa (\hat\Omega)$ which is the unit ball,
whenever $q \in \BB^n (p;\frac{r_3}{2k(m)})$ and $m$ is sufficiently large.
Therefore there exist an integer $m_2 > m_1$ such that
$(L\circ \kappa)^{-1}(\BB^n(0; 1-\epsilon)) \subset\subset \cH_q$ whenever $q \in \BB^n (p;r_3/m_2)$.

As in (\ref{T19}), a direct computation yields
\begin{equation} \label{T23}
\kappa(\cF_\delta) = \{z \in \CC^n \colon |z_1|^2 + \frac{\delta}2 |z_1-\bar z_1| + (1+\delta)H_{b(q)} (z') < 1\}.
\end{equation}
Now, consider the set $L\circ \kappa \circ \Lambda_\lambda (\cF_\delta)$ for each
$\delta > 0$. (Recall that  $\Lambda_\lambda (\cF_\delta)=\cF_\delta$ as remarked
in the line below (\ref{T20}).)  These domains increase monotonically as $\delta \searrow 0$
(since $\cF_\delta$'s do) in such a way that
the union $\bigcup_{0<\delta <\delta_0}  L \circ \kappa \circ (\cF_\delta)$ becomes arbitrarily close to $\BB^n (0; 1)$ as $m$ is sufficiently large.
\begin{figure}[h]
\centering
\includegraphics[height=2in, width=2.70in]{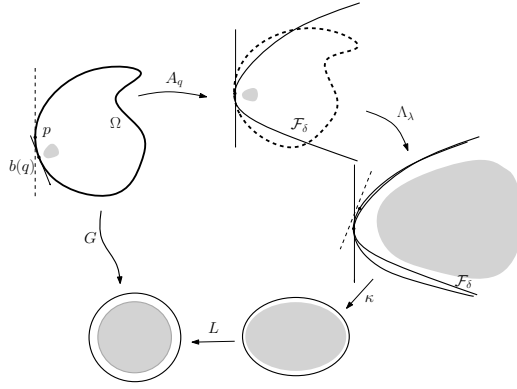}
\caption{$\BB^n (0;1-\epsilon) \subset G(\Omega)$}
\end{figure}
Consequently there exists a constant $\delta>0$ such that $\BB^n(0;1-\epsilon) \subset\subset L \circ \kappa \circ (\cF_\delta)$.  Moreover there is an intger $m_3 > m_2$ such that
\begin{equation} \label{T24}
\Lambda_\lambda^{-1} \big(\kappa^{-1}\circ L^{-1}(\BB^n(0; 1-\epsilon)\big) \subset \BB^n (0;r_3/k(m_1)),
\end{equation}
as $\Lambda_\lambda^{-1}$ scales down the compact subsets (since $\lambda < r_3/m_2$, sufficiently small) to a small set near the origin.
Hence, we have
$$
\Lambda_\lambda^{-1} \big(\kappa^{-1}\circ L^{-1}(\BB^n(0; 1-\epsilon)\big) \subset \cF_\delta \cap \BB^n (0;r_3/k(m_1)) \subset \Omega.
$$
Consequently,
\begin{eqnarray}
\BB^n (0; 1-\epsilon)
& \subset & L \circ \kappa \circ \Lambda_\lambda
(\cF_\delta \cap \BB^n (0;r_3/k(m_1)) ) \nonumber \\
& \subset & L \circ \kappa \circ \Lambda_\lambda (A_q (\Omega)) \label{T25} \\
& =  &G(\Omega), \nonumber
\end{eqnarray}
as long as $q \in \BB^n (p; \frac{r_3}{2k(m_3)})$.
\medskip


Now we show that $G(\Omega) \subset \BB^n (0; 1+\epsilon)$.  Consider
$$
\Omega' := \Omega - \BB^n (p, r_2).
$$
Notice that there exists an integer $\ell >>1$ such that
\begin{equation} \label{T26}
A_q (\Omega') \subset A_q(\Omega) - \BB^n (0; r_2/\ell) \subset \cE - \BB^n (0; r_2/\ell).
\end{equation}

Now, there exists an integer $m_4>3$ such that, if $m>m_4$ and $q \in \BB^n (p,\frac{r_3}{2k(m)})$, then
$$
\Lambda_\lambda (\cE - \BB^n (0; r_2/k))  \subset \{z \in \cE \colon \re z_1 >
\frac{r_2}{r_3}\cdot \frac{m_4}{\ell} \}.
$$
This implies that there exists $m_4$ such that
$$
G(\Omega') \subset L \circ \kappa (\{z \in \cE \colon \re z_1 >
\frac{r_2}{r_3}\cdot \frac{m_4}{\ell} \}) \subset (\BB^n (-{\bf 1}; \rho(m_4)))
$$
for some $\rho(m)$ which approaches zero as $m$ tends to infinity; a direct computation with the Cayley transform and the choice of $L$ (cf.\ (\ref{T14})) verify this immediately.  Therefore, choosing $m_4$ sufficiently large, we arrive at
\begin{equation} \label{T27}
G(\Omega') \subset \BB^n (-{\bf 1}; \epsilon).
\end{equation}
\begin{figure}[h]
\centering
\includegraphics[height=2in, width=3.40in]{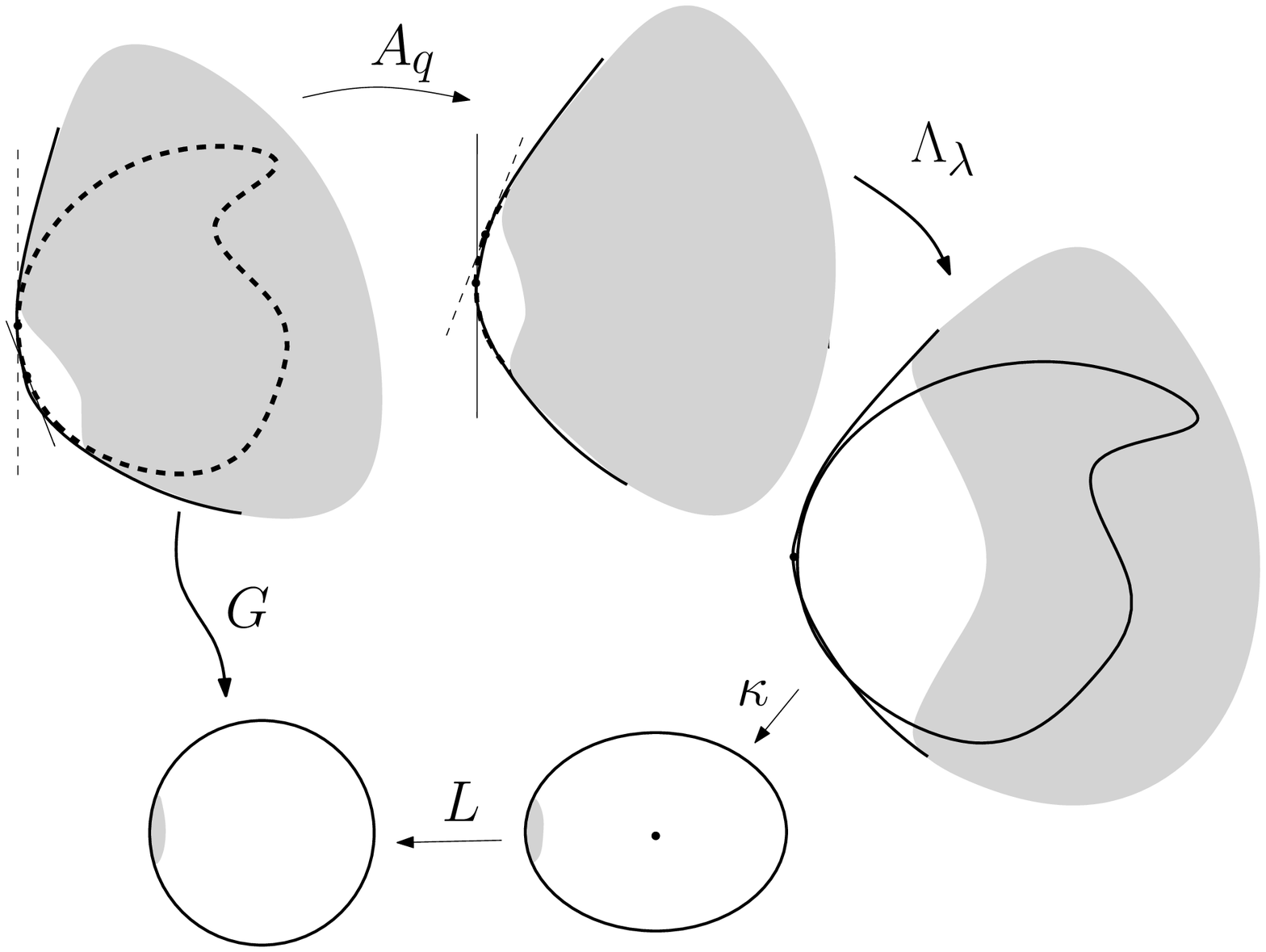}
\caption{$G(\Omega') \subset \BB^n (-{\bf 1};\epsilon)$}
\end{figure}
For the $\epsilon$ given above, there exists $\delta $ such that
\begin{equation} \label{T28}
L\circ \kappa (\cG_\delta) \subset \BB^n (0;1+\epsilon).
\end{equation}
Fix this $\delta$.  Then, recall how the auxiliary domain $\cG_\delta$ was defined
in (\ref{T16}).  Given any $\delta > 0$, according to (\ref{T4})--(\ref{T6}),
there exists $\rho >0$ such that
$$
A_q (\Omega) \cap \BB^n(0;\rho) \subset \cG_\delta.
$$
\begin{figure}[h]
\centering
\includegraphics[height=2in, width=3.40in]{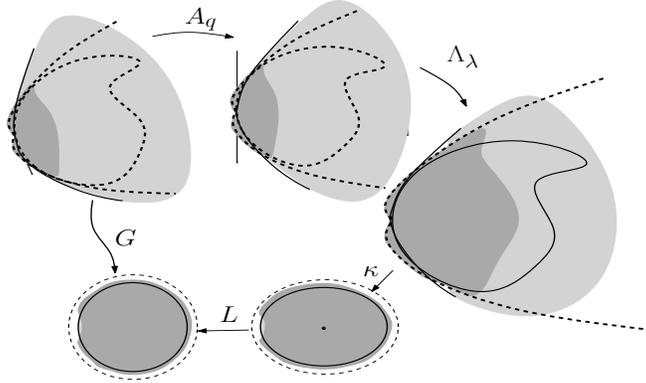}
\caption{$G(\Omega) \subset \BB^n (0; 1+\epsilon)$}
\end{figure}
On the other hand, we can go back to (\ref{T26}) and require that $r_2/\ell < \rho/2$.  Then
we have
\begin{equation} \label{T29}
A_q (\Omega) \cap \BB^n(0;2r_2/\ell) \subset \cG_\delta.
\end{equation}
Since there exists an integer $m_5>0 $ such that
$A_q (\BB^n (p; r_2/\ell) \subset \BB^n (0; 2 r_2/\ell)$, we have that
$$
G(\Omega - \Omega') \subset L \circ \kappa \circ \Lambda_\lambda
\big( A_q(\Omega) \cap \BB^n (0; 2 r_2/\ell) \big).
$$
This implies
\begin{eqnarray}
G(\Omega - \Omega')
& \subset & L \circ \kappa \circ \Lambda_\lambda
\big( A_q(\Omega) \cap \BB^n (0; 2 r_2/\ell) \big) \nonumber \\
& \subset & L \circ \kappa \circ \Lambda_\lambda(\cG_\delta)
\qquad\qquad\qquad\qquad \text{by (\ref{T29})} \label{T30}\\
& \subset & L \circ \kappa (\cG_\delta)
\qquad\qquad\text{by the sentence following (\ref{T20})}\nonumber\\
& \subset & \BB^n (0;1+\epsilon).\nonumber
\end{eqnarray}

By (\ref{T27}) and (\ref{T30}) we have that
$$
G(\Omega) \subset \BB^n (0; 1+\epsilon).
$$
This completes the proofs of Claim and Theorem \ref{thm-2g}.  \hfill $\Box$
\medskip

\section{Remarks}

In this final section we present several remarks.

\subsection{On the spherically-extreme points}
Pertaining to Question \ref{q1}, one of the naturally rising question would be whether
one may re-embed (the closure of) the bounded strongly pseudoconvex domain so that
the pre-selected boundary point becomes spherically extreme.  Recent paper
by Diederich-Fornaess-Wold \cite{DFW} says that the answer to this question is
affirmative.  Owing to this new result, Theorem \ref{thm-2g} now implies the following

\begin{theorem} If $\Omega$ is a bounded domain in $\CC^n$ with a $\cC^2$-smooth strongly pseudoconvex boundary, then $\lim_{\Omega\ni z \to \partial\Omega} \sigma_\Omega (z) = 1$.
\end{theorem}

On the other hand, a more ambitious try may be that one would like to re-embed
the domain using the
automorphisms of $\CC^n$ to achieve the same goal. But this cannot work.  Here is a
counterexample to such a try:

\begin{example} \rm Consider the domain $U$ which is the open $1/10$-
tubular neighborhood  of the circle $S:= \{(e^{it}, 0) \in \CC^2 \colon t \in \RR \}$. This
domain is strongly pseudoconvex.  Let $p = (9/10, 0)$.  Clearly $p \in \partial U$. If
there were $\psi \in \Aut(\CC^2)$ that makes $\psi(p)$ sperically-extreme for
$\psi(U)$, then consider the analytic disc $\Sigma := \psi(\Delta)$ where $\Delta :=
\{(z,0) \colon |z|\le 1\})$. Since $\Delta$ crosses $\partial U$ transversally at $\psi(p)$,
$\Sigma$ crosses the sphere envelope at $\psi(p)$ and extends to the exterior of the
sphere.  On the other hand the boundary of $\Sigma$ remains inside $\psi(U)$ and hence
inside the sphere. Now let the sphere expand radially from its center, and let it stop at the
radius beyond which cannot have intersection with the holomorphic disc $\Sigma$. Then
the sphere is tangent to a point to $\Sigma$ at an interior point keeping the whole disc
inside the sphere.  The maximum principle now implies that $\Sigma$ should be entirely
on the sphere.  But the boundary of $\Sigma$ is strictly inside the sphere, which is a
contradiction. This implies that $p$ cannot be made spherically-extreme via any
re-embedding by an automorphism of $\CC^n$.
\end{example}
\smallskip

\it Acknowledgement: \rm This example was obtained after a valuable discussion between
the first named author and Josip Globevnik.  The first named author would like to
express his thanks to Josip Globevnik for pointing out such possibility.

\subsection{On the exhaustion theorem by Fridman-Ma}

The main theorem by Buma Fridman and Daowei Ma in \cite{FridMa} had obtained the conclusion
of Theorem \ref{thm-2g} in the sepcial case $\Omega \ni q \to p$ {\it trasversely} to the
boundary $\partial \Omega$.  However, that is not sufficient to prove Theorem \ref{thm-2g}; it
is indeed necessary to consider all possible sequences approaching the boundary.
In \cite{FridMa} they need not consider the point sequences approaching the
boundary tangentially, as their interest
was only on the holomorphic exhaustion of the ball by the biholomorphic images of a
bounded strongly pseudoconvex domain.  On the other hand, our proof of Theorem \ref{thm-2g}
gives a proof to their theorem as well; one only need to use $(1+\epsilon)^{-1} G(z)$ instead of
$G$.  [Recall that $G$ depends upon $q$.  Letting $q$ converge to $p$ and $\epsilon$
tend to zero, one gets a sequence of maps that exhausts the unit ball holomorphically.]

\subsection{Plane domain cases}

For domains in $\CC$, several theorems have been obtained by F. Deng, Q. Guan and L.
Zhang in \cite{DGZ}.  Theorem \ref{thm-2g} obviously includes many of those
results, as every boundary point of a plain domain with $\cC^2$ smooth boundary is
spherically-extreme.


\begin{thebibliography}{99}

\bibitem{DGZ} Deng, F.; Guan, Q.; Zhang, L.: On some properties of squeezing
functions of bounded domains, {\it Pacific J. Math.}, 257, no. 2,
(2012), 319--342.

\bibitem{DFW} Diederich, K.; Fornaess, J. E.; Wold, E. F.:
Exposing points on the boundary of a strictly pseudoconvex or a locally convexifiable domain of finite 1-type, arXiv:1303.1976.

\bibitem{FridMa} Fridman, Buma and Ma, Daowei: On exhaustion of domains.
{\it Indiana Univ.\ Math. J.} 44 (1995), no. 2, 385--395.

\bibitem{K} Kim, Kang-Tae: Asymptotic behavior of the curvature of the Bergman
metric of the thin domains, {\it Pacific J. Math.} 155 (1992), no. 1, 99--110.

\bibitem{Klembeck} Klembeck, Paul: K\"ahler metrics of negative curvature, the Bergmann
metric near the boundary, and the Kobayashi metric on smooth bounded
strictly pseudoconvex sets, {\it Indiana Univ.\ Math.\ J.} 27
(1978), no. 2, 275--282.

\bibitem{Lee} Lee, Sunhong: Asymptotic behavior of the Kobayashi metric on certain
infinite-type pseudoconvex domains in $\CC^2$, {\it J. Math.\ Anal.\
Appl.} 256 (2001), no. 1, 190--215.

\bibitem{LSY-1} Liu, Kefeng; Sun, Xiaofeng; Yau, Shing-Tung: Canonical metrics on the
moduli space of Riemann surfaces, I. {\it J. Differential Geom.} 68
(2004), no. 3, 571--637.

\bibitem{LSY-2} Liu, Kefeng; Sun, Xiaofeng; Yau, Shing-Tung: Canonical metrics on the
moduli space of Riemann surfaces, II. {\it J. Differential Geom.} 69
(2005), no. 1, 163--216.

\bibitem{Pinchuk} Pinchuk, Sergey: The scaling method and holomorphic mappings.
Several complex variables and complex geometry, Part 1 (Santa Cruz,
CA, 1989), 151--161, {\it Proc.\ Sympos.\ Pure Math.} 52, Part 1,
Amer.\ Math.\ Soc., Providence, RI, 1991.

\bibitem{Wong} Wong, Bun: Characterization of the unit ball in Cn by its
automorphism group. {\it Invent.\ Math.} 41 (1977), no. 3, 253--257.

\bibitem{Y} Yeung, Sai-Kee: Geometry of domains with the uniform squeezing
property, {\it Adv. Math.} 221 (2009), no. 2, 547--569.



\end{thebibliography}
\end{document}